\documentclass[12pt,a4paper]{amsart}
\usepackage[usenames]{color}
\usepackage{times}
\usepackage{amsfonts}
\usepackage{amsthm}
\usepackage{amsmath}
\usepackage{psfrag}

\usepackage{times}
\usepackage{latexsym,color}
\usepackage{amssymb}
\usepackage{graphicx}

\usepackage{epsfig}

\usepackage[mathscr]{euscript}
\usepackage{tikz}
\usepackage[all]{xy}
 \usepackage[hang]{caption}
 \usetikzlibrary{snakes}
\usetikzlibrary{arrows}

\advance \topmargin by -\headheight
\advance \topmargin by -\headsep
\evensidemargin \oddsidemargin
\newcommand{\GZ}{{\mathbb Z}}

\usepackage{lineno}

\newtheorem{theorem}{Theorem}
\newtheorem{lemma}[theorem]{Lemma}
\newtheorem{corollary}[theorem]{Corollary}
\newtheorem{proposition}[theorem]{Proposition}
\newtheorem{question}{Question}
\newtheorem{example}{Example}
\newtheorem{definition}{Definition}

\newtheorem{conj}{Conjecture}

\begin{document}

\title{Clustering, order conditions, and  languages of interval exchanges}

\author[S. Ferenczi]{S\'ebastien Ferenczi}
\address{Aix Marseille Universit\'e, CNRS, Centrale Marseille, Institut de Math\' ematiques de Marseille, I2M - UMR 7373\\13453 Marseille, France.}
\email{ssferenczi@gmail.com}

\author[L.Q. Zamboni]{Luca Q. Zamboni}
\address{Institut Camille Jordan\\
Universit\'e Claude Bernard Lyon 1\\
43 boulevard du 11 novembre 1918\\
F69622 Villeurbanne Cedex
(France)}
\email{zamboni@math.univ-lyon1.fr}

 \date{March 15,  2026}


\begin{abstract}
We investigate various connections between the clustering for the Burrows-Wheeler transform, a lossless algorithm used in data compression, and languages of interval exchange transformations. We show that a primitive word $u$ clusters for a pair of orders $(<_D,<_A)$ if and only if $u$ is a return word in the natural coding of a generalised interval exchange transformation with departure and arrival orders $(<_D,<_A)$. This answers a question of M. Lapointe on the perfect clustering of return words for a symmetric standard interval exchange transformation. We show that if $T$ is symmetric, then all natural codings are palindromically rich languages, and the orders of the induced transformation on a cylinder $[w]$ equal the original departure/arrival orders $(<_D,<_A)$ if and only if the shortest bispecial word containing $w$ is a palindrome. We also investigate language related features distinguishing between standard and generalised interval exchange transformations.
 \end{abstract}
\maketitle

In this paper we investigate various connections between  {\it clustering} for the {\it Burrows-Wheeler transform}, a lossless compression algorithm used in various applications including image and data compression and of purely word combinatorial nature, and an important class of one dimensional dynamical systems, the {\it interval exchange transformations}. 

Given a word $u$ over an ordered  alphabet $(\mathcal A, <),$  the cyclic conjugates of $u$ are arranged in increasing lexicographic order thereby forming a rectangular array of words with the lexicographically  smallest at the  top and the  largest at the bottom. Thus the first column of this array reading from top to bottom is of the form $a_1^{n_1}a_2^{n_2}\cdots a_k^{n_k}$ with $n_i\geq 1$  and  $a_1<a_2 <\cdots < a_k.$  The Burrows-Wheeler transform of $u$ is defined as the last column of this array of words. This process defines a  lossless compression algorithm in that the original word $u$ may be completely recovered up to cyclic conjugacy from its Burrows-Wheeler transform.    A word $u$ is said to {\it cluster} if its Burrows-Wheeler transform  is of the form $b_1^{m_1}b_2^{m_2}\cdots b_k^{m_k}$ with $m_i\geq 1$ and  distinct letters $b_i.$     In other words, a word $u$ clusters if and only if the last column of the lexicographic array of $u$ is obtained by some permutation of the letter blocks in the first column. Thus each clustering word $u$ over an ordered alphabet $\mathcal A$ defines a permutation of $\mathcal A$  or equivalently a second order  on $\mathcal A$ determined by the order of the letter blocks in the last column.  Alternatively, as shall be our convention in the present paper, we  start with a fixed pair of orders
$(<_D, <_A)$ on the alphabet $\mathcal A $, called the {\it departure} and  {\it arrival} orders respectively, and we say that a word $u$ {\it clusters} for the pair  $(<_D, <_A)$ if by arranging the conjugates of $u$ in lexicographic order with respect to the departure order $<_D,$ the letters in the last column are arranged according to the arrival order $<_A.$ \\

The natural coding of a $k$-interval exchange transformation $T$  for  a pair  $(<_D, <_A)$  is given by a mapping $T$ on an interval $I$ partitioned into $k$ subintervals labelled by the letters in $\mathcal A$ in increasing order with respect to $<_D$   while their images  under $T,$ which also partition $I,$ are labelled in increasing order with respect to  $<_A.$ Depending on the nature of the interval exchange transformation one imposes different conditions on the restriction of $T$ to each of the $k$ labelled subintervals: in a {\it standard interval exchange transformation} each restriction is  a translation, in an {\it affine  interval exchange transformation} each  restriction is  an affine mapping and in a {\it generalised interval exchange transformation} each restriction is  assumed only to be continuous and (strictly) increasing. 

A common feature linking clustering words  and natural codings of interval exchange transformations is the following  {\it order condition} on a language $L$,  defined in terms of a pair of orders  $(<_D, <_A)$ over the underlying alphabet: whenever $axb, cxd \in L$ with $a, b,c,d \in \mathcal A,$  $x\in \mathcal A^*$ (possibly empty), and $a\neq c$ and $b\neq d,$ then $a<_A c$ if and only if $b<_D d.$ \\   

Our starting point in the present paper concerns  a question of M. Lapointe \cite{lap} in which she asks whether each return word in a natural coding of a {\it symmetric} standard interval exchange transformation is perfectly clustering. Symmetric here means that  the defining orders $(<_D, <_A)$ are symmetric in the sense that for each pair of letters $a,b\in \mathcal A$ we have $a<_D b$ if and only if $b<_A a,$ and similarly perfectly clustering means clustering for some  pair  of symmetric orders $(<_D, <_A).$ Lapointe's question was answered in part by F. Dolce and C.B. Hughes \cite{d}, where it is shown that if $u$ is a return word in a regular interval exchange transformation $T$ (not necessarily symmetric), then $u$ clusters for some pair of orders. Regular here refers to a standard interval exchange transformation satisfying Keane's i.d.o.c. or infinite disjoint orbit condition \cite{kea}.  However, the clustering orders remained unknown, so while they could show that each return word in a symmetric regular interval exchange transformation clusters, it remained an open question whether it perfectly clusters.  Finally, Dolce and Hughes ask in \cite{d}  whether return words in a standard interval exchange transformation for a fixed pair of orders cluster for the same pair of orders.  \\

In the present paper we settle these questions by showing that for a fixed pair of orders $(<_D, <_A)$ on an alphabet $\mathcal A,$ a primitive word $u$ over $\mathcal A$ clusters for the pair $(<_D,<_A)$ if and only if $u$ is a return word in a natural coding of a generalised interval exchange transformation (not necessarily standard or regular) defined by the pair $(<_D,<_A).$  Thus in particular each return word in a symmetric (standard, affine or generalised) interval exchange transformation is perfectly clustering. Our proof is short and combines two earlier results of the authors: the first result (Theorem 1 in  \cite{fz5})  states that a word $u$ over $\mathcal A$ clusters for a pair $(<_D, <_A)$ if and only if the language $L_u$ consisting of all factors of $uuu \cdots$ satisfies the order condition for the pair of orders $(<_D,<_A).$ The second result, obtained in collaboration with P. Hubert  in \cite{fhuz}, states that an extendable language  $L$ is a natural coding of a  generalised interval exchange transformation with departure and arrival orders  $(<_D,<_A)$ if and only if $L$ satisfies the order condition for the pair of orders $(<_D,<_A).$ 
 
 We then extend this result to concatenations of return words. More precisely, given a pair of orders   $(<_D,<_A)$ on an alphabet $\mathcal A$ and a language $L$ over $\mathcal A$ satisfying the order condition for $(<_D,<_A),$  we consider the set $R_w$ of all  return words in $L$ to some fixed word $w\in L.$ We then define a new alphabet $\mathcal A'$ in which each letter in $\mathcal A'$ corresponds to some return word  $u$ in $R_w.$ More precisely we fix a bijection $\phi: \mathcal A' \rightarrow R_w$ and extend $\phi$ to a morphism from the free monoid generated by $\mathcal A'$ to the free monoid generated by $\mathcal A.$  We also extend the orders $<_D$ and $<_A$ to primitive words $u$ and $v$ over $\mathcal A$ as follows: we write $u<_D v$  if $uuu\cdots $ is lexicographically smaller than $vvv\cdots$ and $u<_A v$ if $\cdots uuu$ is reverse lexicographically (meaning reading from right to left) smaller than $\cdots vvv.$ This allows us to define a pair of induced orders $(<_{Dw}, <_{Aw})$  on $\mathcal A'$ given  by $a'<_{Dw} b' $ whenever $\phi(a') <_D \phi(b')$ and similarly  $a'<_{Aw} b' $ whenever $\phi(a') <_A \phi(b').$  When  the induced interval exchange transformation $T_{[w]}$ on the cylinder $[w]$ exists, these are precisely its departure and arrival orders: the cylinder $[w]$ is partitioned according to cylinders of  words in $R_w,$ these cylinders are labelled by letters in $\mathcal A'$ in increasing order with respect to $<_{Dw},$ and $T_{[w]}$ permutes these cylinders according to the induced arrival order $<_{Aw}.$ In Theorem \ref{t2} we prove that a word $u'$ over the alphabet $\mathcal A'$ clusters for $(<_{Dw}, <_{Aw})$ if and only if $\phi(u')$ clusters for $(<_D, <_A).$  In particular, we recover our earlier result  that each return word $u\in R_w$ clusters for $(<_D, <_A):$ writing $u=\phi(a')$ for some letter $a'\in \mathcal A',$  since $a'$ clusters for $(<_{Dw}, <_{Aw})$ (each letter trivially clusters for any pair of orders),  we get that $u$ clusters for $(<_D, <_A).$ 
 We further prove that if $T$ is minimal, aperiodic and the return word $u\in R_w$ is sufficiently long, then $u$ clusters for no pair of orders other than $(<_D, <_A).$ Finally we  show that our results here do not extend to interval exchange transformations with {\it flips}. \\

If we now fix a natural coding of an interval exchange transformation $T$ for a pair of orders $(<_D,<_A),$ we may ask for which words $w$ are the derived orders $(<_{Dw}, <_{Aw})$ on the induced transformation $T_{[w]}$ equal to the original pair $(<_D,<_A).$ 
In the important case of standard symmetric interval exchange transformations, we show that the induced pair of orders is equal to the symmetric pair  $(<_D,<_A)$ if and only if the shortest bispecial extension of $w$ (i.e., the shortest bispecial word containing $w)$ is a palindrome.
Our proof uses the full combinatorial theory of {\it rich languages}, those languages which contain, in some sense, an abundance of palindromes. Namely, a language $L$ is rich if each complete return word to a palindrome is a palindrome. We prove that if $L$ is a language satisfying the order  condition for some pair of symmetric orders $(<_D, <_A)$ then $L$ is a rich language.  In particular each natural coding of a symmetric interval exchange transformation is a rich language. We  give also some estimates on the {\it palindromic complexity} of  $L,$ i.e., the number of palindromes in $L$  of each given length $n.$ \\
  
 The notion of clustering gives new insights on a difficult and much-studied question, namely to understand the difference in behaviour between standard and generalised  interval exchange transformations, see \cite{fhuz} for references. The clustering of a primitive word $w$ is linked to an order condition on the infinite language $L_w$ (or equivalently on the finite language consisting of all  factors of $ww)$ and it follows from \cite{fz4} and  \cite{fz5} that   each word in this language is produced by a natural coding of a standard interval exchange transformation if and only if it is produced by a natural coding of a generalised interval exchange transformation. In order to
 distinguish between natural codings of standard and generalised interval exchange transformations, we are led to consider  the order condition only on the set of factors of a word $w$, and not $ww,$ the two clearly not being equivalent (e.g., $w=abba).$    We show that for each finite set of words $W,$ if the (non extendable) language consisting of all factors of words in $W$ satisfies the order condition  
 for some pair $(<_D,<_A),$ then $W$ is contained in the natural coding of  a generalised interval exchange transformation, but not necessarily that of a standard one. If $(<_D,<_A)$ are symmetric orders, then, using the theory of  self-dual induction from \cite{fz1}, we are able to determine which sets $W$ satisfying the order condition for $(<_D,<_A)$ are contained in the natural coding of a standard interval exchange transformation with these orders: sets $W$ made with a single word are among those, while this result is no longer true when $|W|\geq 2.$

 \section{Languages}
 Let $\mathcal A$ be a finite set called the {\em  alphabet}, its elements being {\em letters}.
 A {\em word} $w$ of {\em length} $n=|w|$ is $w_1w_2\cdots w_n$, with $w_i \in {\mathcal A}$. The set of words on $\mathcal A$ is denoted by $\mathcal A^{\star}$. The {\em concatenation} of two words $w$ and $w'$ is denoted by $ww'$.\\
 A word is  {\em primitive} if it is not a power of a shorter word.\\
 The {\em reverse} of a word $w=w_1 \cdots w_n$
 is the word $\bar w=w_n\cdots w_1$. If $w=\bar w$, $w$ is a {\it palindrome}.\\
 
By a   language $L$ over $\mathcal A$ we mean  a {\it factorial  language}:
a collection of sets $(L_n)_{n\geq 0}$ where the only element of $L_0$ is the {\em empty word} $\varepsilon$, and where each $L_n$ for $n\geq 1$ consists of words of length $n$,  such that
each $v\in L_{n+1}$ can be written in the form $v=au=u'b$ with $a,b\in \mathcal A$ and $u,u'\in
L_n.$ Except when otherwise mentioned, all our languages are tacitly assumed to be {\it extendable}: for each $v\in L_n$ there exists $a,b\in \mathcal A$ with $av,vb\in L_{n+1}$.\\
A word $v=v_1\cdots v_r$  {\em occurs} at index $i$ in a word $w=w_1\cdots w_s$  if $v_1=w_i$, \ldots ,$v_r=w_{i+r-1}$, we say also that $w$ contains $v$ and $v$ is a {\em factor} of $w$.
The  occurrences of a word $v$ in a word $w$ can be as a {\em prefix}, as a {\em suffix}, or {\em internal} occurrences.\\

A language $L$ is {\em closed under reversal} if $w\in L \Leftrightarrow \bar w \in L$.\\
A language $L$  is {\em aperiodic} if for  all nonempty words $w$ in $L$, there exists $n$ such that  $w^n$ is not in $L$.\\
 A language $L$ is
{\em recurrent} if for each $w\in L$ there exists  a nonempty $v\in L$, such that $wv$ ends with $w$.\\
A language $L$ is {\em uniformly recurrent} if for every word $w$ in $L$, there exists a constant $K$ such that $w$ occurs in every word in $L$ of length at least $K$.\\

 A word $w$ in $L$ is {\em  right special} (resp. {\em  left special}) if it has more than one {\em right extension} $wx$ (resp. {\em left extension} $xw$) in $L,$ with $x$ in $\mathcal A.$  If $w$ is both right special and
left special, then $w$ is {\em  bispecial}. If $\# L_1>1$, the empty word $\varepsilon$
 is bispecial.
 To {\em resolve} a bispecial word $w$ is to find all its  {\em bilateral  extensions}  $xwy$ in $L$ with $x$ in $\mathcal A,$  $y$ in $\mathcal A.$\\
 A bispecial word $w$ in $L$ is a {\em strong  bispecial} if
 $$\#\{aw \in L, a\in \mathcal A, b\in \mathcal A\} > \#\{a \in L, a\in \mathcal A\} +\#\{wb \in L, b\in \mathcal A\} -1.$$
 
 A {\it suffix return word} to a word $w$ in a language $L$ is any word $v$ in $L$ such that $wv$ has exactly two occurrences of $w$, one as a prefix and one as a suffix. A {\it prefix return word} to a word $w$ in a language $L$ is any word $v$ in $L$ such that $vw$ has exactly two occurrences of $w$, one as a prefix and one as a suffix. A {\it complete return word} of $w$ is $wv$ for a suffix return word $v$ or $vw$ for a prefix return word $v$.\\

For a word $u$, we denote by $u^{\omega+}$ the one-sided infinite word $uuu\cdots$, by $u^{\omega-}$ the one-sided infinite word $\cdots uuu$, by $u^{\omega}$ the two-sided infinite word $\cdots uuu\cdots$, and by $L_u$ the language consisting of all the factors of  $u^{\omega+}$, or equivalently all the factors of  $u^{\omega-}$, or equivalently all the factors of  $u^{\omega}$.

If $W$ is a set a words on $\mathcal A$, $F(W)$ is the non extendable language made with all the factors of all the  words in $W$. If $W$ has a single element $w$, we denote $F(W)$ by $F(w)$.

\section{Interval exchange transformations}
\subsection{The transformations}
All intervals we consider are closed on the left and open on the right.

\begin{definition}\label{dit} Let $\mathcal A$ be a finite alphabet. A {\em generalised interval exchange transformation}  is  a map $T$ defined on an interval $I$, partitioned into intervals  $I_e$, $e \in  \mathcal A$, continuous and (strictly) increasing on each $I_e$,  and 
such that the $TI_e$, $e\in \mathcal A$, are  intervals partitioning $I$. 

The $I_e$, indexed in $\mathcal A$, are called the {\em defining intervals} of $T$. The $I_e$ will always be denoted, from left to right, by $[\gamma_i, \gamma_{i+1})$, $0\leq i\leq k-1$, and the  $TI_e$ will always be denoted, from left to right, by $[\beta_i, \beta_{i+1})$, $0\leq i\leq k-1$.

If the restriction of $T$ to each $I_e$ is an affine map, $T$ is an {\em affine} interval exchange transformation.

If the restriction of $T$ to each $I_e$ is an affine map of slope $1$, $T$ is a {\em standard} interval exchange transformation.

$T$ is {\em aperiodic} if $T^nx\neq x$ for all $x$ and $n$.  $T$ is {\em minimal} if all its orbits are dense. 
\end{definition}

In this definition, the defining intervals $I_e$ are not necessarily the intervals of continuity of $T$, as it may happen that $I_e$ and $I_f$ are adjacent and sent to adjacent intervals in the correct order, thus $T$ may be extended to a continuous map on their union. Thus in the present paper an interval exchange transformation $T$ is always supposed to be given {\em together with its defining intervals} as keeping the same $T$ but changing the defining intervals would change the coding.

  \begin{definition}
  Let $T$ be a generalised  interval exchange on $k$ intervals, defined on
 an interval $I$, $E$ a subset of $I$. For $x$ in $E$, let $n(x)$ be the smallest $n>0$ such that $T^nx$ is in $E$;  if $n(x)$ is finite for 
 every $x\in E$, the {\em induced map} $T_{E}$ of $T$ on $E$ is the transformation on $E$ defined by $T_Ex=T^{n(x)}x$.
 \end{definition}

 \begin{definition}\label{diord} A generalised interval exchange transformation defines an (oriented) pair of  orders $(<_D, <_A)$   on the alphabet  $\mathcal  A$, respectively called the departure and arrival order, by:
\begin{itemize}
\item $e<_D f$ whenever the interval $I_e$ is strictly to the left of the interval $I_f$,
\item $e<_A f$ whenever the interval $TI_e$ is strictly to the left of the interval $TI_f$.
\end{itemize}
A pair of orders $(<_D, <_A)$ is {\em symmetric} if $<_A$ is the reverse of $<_D$,  i.e.   $a<_A b$ if and only if $b<_D a$ for all distinct letters $a, b\in \mathcal A.$
$T$ is {\em symmetric} if  its pair of orders is symmetric.
\end{definition}
 The orders we consider on alphabets are always strict and total.\\
The orders defined by $T$ correspond to the two {\it permutations} used by Kerckhoff \cite{ker} to define
 standard interval exchange transformations: the unit interval is partitioned into semi-open intervals
which are numbered from $1$ to $k$, ordered according to a permutation  $\pi_0$ and then rearranged according to another permutation  $\pi_1$; then  $x<_A y$ if $\pi_1^{-1} x< \pi_1^{-1}  y$, $x<_D y$ if $\pi_0^{-1} x< \pi_0^{-1}  y$; in more classical definitions, there is only one permutation $\pi$, which corresponds to $\pi_1$ while $\pi_0=Id$; note that in some papers the rearranging is by $\pi_0^{-1}$ and $\pi_1^{-1}$.

\subsection{The codings}
\begin{definition}\label{sy} For a generalised 
interval exchange transformation $T$, coded by $\mathcal A$, its {\em natural coding} is  the  language
 $L(T)$ generated by the {\em trajectories} of all points in $I$,  where the trajectory of a point $x$ in $I$ is the
sequence
$(x_n, n\in\GZ) \in \mathcal A^{\GZ}$ where $x_{n}=e$ if ${T}^nx$ falls into the defining interval 
$I_e$ of $T$ $e\in \mathcal A$. For a word $w=w_0\cdots w_{m-1}$ in $L(T)$, the {\em cylinder} $[w]$ is the  set $\{x\in I; x_0=w_0, ..., x_{m-1}=w_{m-1}\}$.

We say a finite set of words $W$ is {\em produced} by the   generalised 
interval exchange transformation $T$ if it is included in $L(T)$.
\end{definition}

We recall an old result, which can be traced back to \cite{kat}, though it is proved only for standard $k$-interval exchange transformations, but the reasoning is the same for generalised $k$-interval exchange transformations: for a cylinder $[w]$, if  $T_{[w]}$ exists, it is a generalised interval exchange on at most $k$ intervals. If there is no connection in $L(T)$, $T_{[w]}$ is a generalised interval exchange on exactly  $k$ intervals. If $T$  is a standard interval exchange transformation,  $T_{[w]}$ always exists and is also a standard interval exchange transformation.\\

The natural coding of an interval exchange transformation is its coding by the defining intervals. However, in parts of Section 5, for affine interval exchange transformations, we shall need a slightly different notion, where several adjacent defining intervals may be coded by the same letter.

\begin{definition}\label{gc} A language $L$ is a {\em grouped coding} of an affine interval exchange transformation $T$ if there exist intervals $\tilde I_e$, $e\in\tilde{\mathcal A}$ such that
\begin{itemize}
\item each $\tilde I_e$ is an interval, and a disjoint  union of defining intervals of $T$,
\item $T$ is a continuous monotone map on each $\tilde I_e$,
\item $L$ is the language generated by the trajectories $(x_n, n\in \GZ) \in \tilde{\mathcal A}^{\GZ}$ where $x_{n}=e$ if ${T}^nx$ falls into
$\tilde I_e$, $e\in \tilde{\mathcal A}$, for all points $x$ in $I$.
\end{itemize}
\end{definition}

\subsection{The order condition}
The following condition was first mentioned under that name, for the symmetric case, in \cite{zdl}, where it is related to the Markoff condition \cite{mar}, and \cite{fhuz} in the general case.  It appears also in the {\em ordered dendrics} of \cite{d} and the {\it planar trees} of \cite{b+}.

 \begin{definition}\label{dfordc}
  A language $L$ (extendable or not) satisfies  the {\em order condition} for the pair of orders $(<_D,<_A)$ if whenever $axb, cxd \in L$ with $a, b,c,d \in \mathcal A,$  $x\in \mathcal A^*$ (possibly empty),  $a\neq c$ and $b\neq d,$ then $a<_A c$ if and only if $b<_D d.$ \\   
The order condition is called the  {\it symmetric order condition} if the  pair of orders $(<_D,<_A)$ is symmetric. \end{definition} 
 
 Note that the order condition has to be checked only when the $x$ above  is bispecial, and then we say that $x$ in $L$ satisfies the order condition, for the orders $(<_D,<_A)$.

The order condition is linked to  interval exchange transformations by the following result from \cite{fhuz} (Proposition 8 in the easy direction, and Theorem 19 in the difficult one); in many cases it could be deduced from the earlier \cite{bel}.

\begin{theorem}\label{t0}\cite{fhuz} {\em An (extendable) language  $L$ is the natural coding of a  generalised interval exchange transformation for a  pair of orders $(<_D,<_A)$, if and only if $L$ satisfies the order condition for the pair of orders $(<_D,<_A)$.} \end{theorem}

   \begin{definition}\label{con}
If a language $L$ on an alphabet $\mathcal A$ satisfies the order condition for $(<_D,<_A)$, a bispecial word $w$ has a {\em connection} if there are letters $a<_{A}a'$, consecutive in the order $<_{A}$, letters $b<_{D}b'$, consecutive in the order $<_{D}$, such 
that $awb$ and  $a'wb'$ are in $L$, and neither $awb'$ nor $a'wb$ is in $L$.
  \end{definition}
  
If the natural coding of an  interval exchange transformation $T$ has no connection, we say that $T$ has no connection.  By \cite{kea}, if a standard  interval exchange transformation has no connection, it is aperiodic and  minimal and its natural coding  is aperiodic and uniformly recurrent. 

  \section{Clustering}
 The Burrows-Wheeler transform was  defined in \cite{bw}; the image of a word $w$ is another word $B(w)$ which has  in general more repetitions of letters that $w$, and from which $w$ can be recovered up to cyclic conjugacy. Thus in data compression it can be useful to apply $B$ to $w$ before compressing it. This is 
 particularly true when all identical letters in  $B(w)$ are grouped together,  thus the interest of the notion of clustering. 
 
\begin{definition}  Let  $\mathcal A$. be an alphabet, ordered by an order denoted by $<_D$, thus $\mathcal A=\{x_1<_D \ldots <_D x_r\}$. The {\em (cyclic) conjugates} of $w$ are  the words $w_i\cdots w_nw_1\cdots w_{i-1}$, $1\leq i\leq n.$ If $w$ is primitive,  $w$ has precisely $n$ conjugates, and $w_{i,1}\cdots w_{i,n}$ 
denotes the $i$-th  conjugate of $w$ where the $n$ conjugates of $w$ are ordered by ascending lexicographical order.\\
 The {\em  Burrows-Wheeler transform} of $w$ is the word $B(w)=w_{1,n}w_{2,n}\cdots w_{n,n}.$
 In other words, $B(w)$ is obtained from $w$ by first ordering its conjugates in ascending order in a rectangular array, and then reading off the last column.\\
 We say $w$ is {\em clustering for the orders} $(<_D,<_A)$  if $B(w)=(x_{i_1})^{n_{i_1}}\cdots (x_{i_r})^{n_{i_r}}$, where $x_{i_1} <_A \cdots <_A x_{i_r}$ and $n_a$ is the number of occurrences of $a$ in $w$ (we allow some of the $n_a$ to be $0$, thus, given  $w$, there may be several 
 possible pairs of orders). $w$ is {\em perfectly} clustering if it is clustering for a symmetric pair of orders.
 \end{definition}
 
Note that in the usual definitions the order $<_D$ is understated, and denoted just by $<$, while the order $<_A$ is defined by a permutation $\pi$, with $x<_A y$ whenever $\pi^{-1} x< \pi^{-1}  y$.\\

{\it Non-primitive words}. As remarked in \cite{mrs}, the Burrows-Wheeler transform can be extended to a non-primitive word $w_1\cdots w_n$, by ordering its $n$ (non necessarily distinct) cyclic conjugates  by non-strictly increasing lexicographical order and taking the word made by their last letters. Then $B(v^m)$ is deduced from $B(v)$ by replacing each of its letters $x_i$ by $x_i^m$, and $v^m$ is clustering for  $\pi$ if and only if $v$ is clustering for  $\pi$.\\

The clustering is linked to the order condition by the following result, which is a restatement of Theorem 1 of \cite{fz5}, omitting the hypothesis of primitivity as explained in \cite{fz5} just after its proof:

 \begin{theorem}\label{t1} \cite{fz5} A  word $u$ over $\mathcal A$ is clustering for the pair of orders  $(<_D,<_A)$ if and only if  the language $L_u$  satisfies the order condition for  $(<_D,<_A)$. \end{theorem}
 
 The following criterion will be useful.
 
  \begin{lemma}\label{useful} For a word $u$ over $\mathcal A,$ the language $L_u$ satisfies the order condition for a pair of orders $(<_D,<_A)$ if and only if  for all $axb,cxd$ occurring in $uu$ with $a\neq c$ and $b\neq d,$ we have $a<_A c\Leftrightarrow b<_D d.$
\end{lemma}

\begin{proof} By the remark after Definition \ref{dfordc},  all we have to prove is that all bispecials in $L_u$ and their bilateral extensions are factors of $uu$. This is proved for primitive words in Theorem 1 of \cite{fz5}, while, for $u=v^k$, $L_u=L_v$ and $vv$ is a factor of $uu$.  \end{proof}

 We shall now show that return words of  interval exchange transformations cluster.

\begin{proposition}\label{fz7}Let $L$ be a language satisfying the order condition for a pair of orders $(<_D, <_A).$  Let $u\in L$ be a return word to some word $w\in L.$ Then $L_u$ satisfies the order condition for $(<_D,<_A).$
\end{proposition}
\begin{proof} Without loss of generality we may assume that $u$ is a prefix return to a word $w.$ Hence $uw\in L$ begins and ends in $w$ and has no other occurrences of $w.$ To show that $L_u$ verifies the  order condition for $(<_D,<_A)$, by Lemma~\ref{useful} it suffices to show that if $axb, cxd$ are factors of $uu$ with $a\neq c$ and $b\neq d,$ then $a<_A c \Leftrightarrow b<_D d.$ As $axb$ and $cxd$ each occur in $uu,$ it follows that there are at least two internal occurrences of $x$ in $uu.$ But  $uu$ has only one internal occurrence of $w$ and hence $x$ does not contain $w$ as a factor. Therefore $axb$ and $cxd$ are each factors of $uw\in L$ and hence $axb, cxd \in L.$ As $L$ satisfies the order condition for $(<_D,<_A),$ it follows that $a<_A c\Leftrightarrow b<_D d$ as required. \end{proof}

\begin{theorem}\label{simple} Let $(<_D, <_A)$ be a pair of orders on an alphabet $\mathcal A$ and let $u$ be a primitive word over $\mathcal A.$ Then $u$ clusters for the pair $(<_D, <_A)$ if and only if $u$ is a return word in a natural coding of a generalised interval exchange transformation with orders $(<_D, <_A).$ \end{theorem}
 \begin{proof} 
  First assume that $u$ clusters for the pair $(<_D, <_A).$ Then by Theorem~\ref{t1}, the language $L_u$ satisfies the order condition for $(<_D, <_A).$ Moreover,  the language $L_u$ is extendable and $u$ is a return word in $L_u$ since $uu$ begins and ends in $u$ and, as $u$ is primitive, there are no internal occurrences of $u.$ And therefore, by Theorem~\ref{t0}, the word $u$ is a return word in a natural coding of a generalised interval exchange transformation. 
 
 For the converse, suppose $u$ is a return word in the natural coding of a generalised  interval exchange transformation $T$ with pair of orders $(<_D, <_A).$  Then by Theorem~\ref{t0}, the language $L$ of $T$ satisfies the
 order condition for the pair   $(<_D, <_A),$  and hence by Proposition~\ref{fz7} the language  $L_u$ satisfies the order condition for $(<_D,<_A). $ It now follows from Theorem~\ref{t1} that $u$ clusters for the pair $(<_D,<_A)$ as required. \end{proof} 
 
Note that a prefix or suffix return word is always primitive: this was proved in Remark 4.5.10 of \cite{casn}, but is also  a direct consequence of the fact that for any word $w$ $R_w$ is a code, which was proved much earlier in \cite{dhs}; indeed, the ``if" direction of Theorem \ref{simple} does not use the  primitivity. Thus in particular each return word in a symmetric (standard or generalised)  interval exchange transformation perfectly clusters, which  answers the question of Lapointe in \cite{lap}.\\

We next extend Theorem~\ref{simple} to concatenations of return words.

\begin{theorem}\label{t2} Let  $\mathcal A$ be an alphabet with $k$ letters, with a pair of orders $(<_D,<_A)$. Let $L$ be a language on $\mathcal A$ satisfying the order condition for  $(<_D,<_A)$. Let $w$ be a word in $L$. Suppose the set $R_w$ of the suffix (resp. prefix) return words to  $w$  is nonempty; let
$R_w=\{U_1, ... , U_s\}.$ Let $\phi$ be  a morphism from $\mathcal A'=\{a_1, \ldots, a_s\}$ 
to $\mathcal A^{\star}$ given  by $\phi a_i=U_i$. \\
Then there exist two orders on  $\mathcal A'$, which we call the {\em derived orders with respect to $w$}, defined as follows, if the $V_l$ are the complete return words, namely $V_l=w\phi a_l$ if we consider suffix return words,  $V_l=\phi a_lw$ if we consider prefix return words, 
\begin{itemize}
\item $a_i<_{Dw} a_j$ if $(V_i)_t<_D(V_j)_t$ where $t$ is the smallest index on which the words  $V_i$ and $V_j$ differ.
\item $a_i<_ {Aw}a_j$ if $(V_i)_u<_A(V_j)_u$  where $u$ is the largest index on which the words $V_i$ and $V_j$ differ.
\end{itemize}
And for any word $v$ on $\mathcal A'$,  $\phi v$ clusters for $(<_D,<_A)$ on $\mathcal A$  if and only if $v$ clusters for $(<_ {Dw},<_ {Aw})$ on $\mathcal A'$.
\end{theorem}
\begin{proof}
We consider suffix return words. Let $U_{i_1}\ldots ,U_{i_r}$ be any sequence of  return words to $w$. Then  $wU_{i_1}=V'_1w, wU_{i_2}=V'_2w, \ldots, wU_{i_r}=V'_rw$; this describes
all the possible occurrences of $w$ in the word
 $wU_{i_1}\cdots U_{i_r}$, and there are $r+1$ of them. Moreover,  $wU_{i_1}\cdots U_{i_r}$ has $w$ as a suffix, and so does the word
 $U_{i_1}...U_{i_r}$  whenever it is at least as long as $w$.

 When $i\neq j$, $wU_i$ is not a strict prefix of $wU_j$, resp.  $wU_j$ of $wU_i$, otherwise there is a forbidden occurrence of $w$ in $wU_j$, resp. $wU_i.$ Thus the index $t$ defining $<_{Dw}$ exists, the order $<_ {Dw}$ coincides with the lexicographical order defined by $<_D$ for the $U_i^{\omega+}$  
 and $<_D$ is indeed an order.    When $i\neq j$, $wU_i$ is not a strict suffix of $wU_j$, resp.  $wU_j$ of $wU_i$, otherwise there is a forbidden occurrence of $w$ in $wU_j$, resp. $wU_i.$ Thus the index $u$ defining  $<_{Aw}$ exists,  the order $<_ {Aw}$ coincides with the reverse lexicographical order defined by $<_A$ for the $U_i^{\omega-}$  
and $<_A$ is indeed an order.   \\

Suppose first every  bispecial in $L_v$ satisfies the order condition for $(<_{Dw},<_{Aw})$, and let $X$ be a bispecial word in  $L_{\phi v}$, with extensions $xXy$ and $x'Xy'$, $x,x',y,y'$ letters of $\mathcal A$, $x\neq x'$, $y\neq y'$. 
 If $X$ contains no occurrence of $w$, then 
$X$ and its extensions must be contained in some $w\phi a_i$, thus in $L$, thus $X$ satisfies the order condition for $(<_D,<_A)$. 

Otherwise, we look at the occurrences of $w$ in $X$. Let the first one be $X_c.\cdots X_{c+g}$, the last one $X_d\cdots X_{d+g}$, with possibly $c=d$. We can then identify uniquely $X_{c}...X_{d+g}$ as either $w$ or
 $w\phi a_{i_1}...\phi a_{i_r}$.
Then $xXy$ and $x'Xy'$
end respectively with $X_{d+g+1}\cdots X_my$ and $X_{d+g+1}\cdots X_my'$, which, as $y\neq  y'$,  are respective prefixes of $\phi a_j$ and  $\phi a_{j'}$, $j\neq j'$; by definition of the orders, $y<_D y'$ if and only if $a_{j}<_{Dw} a_{j'}$. 
$xXy$ and $x'Xy'$
begin with $xX_{1}\cdots X_{c+g}$ and  $x'X_{1}\cdots X_{c+g}$. which because of the structure of $L_{\phi v}$ are respective suffixes of some $w\phi a_i$ and $w\phi a_{i'}$, and $i\neq i'$ as $x\neq x'$; by definition of the orders,  $x<_A x'$ if and only if $a_i<_{Aw} a_{i'}$. Thus we have proved that $a_{i_1}\cdots a_{i_r}$
 (possibly empty) is  bispecial in $L_v$, and  thus by the hypothesis satisfies the order condition for $(<_{Dw},<_{Aw})$, hence $X$ satisfies the order condition for $(<_D,<_A)$.\\

Suppose now every bispecial in $L_{\phi v}$ satisfies the order condition for $(<_D,<_A)$, and let $X$ be a bispecial word possibly empty) in $L_v$, with extensions $a_iXa_j$ and $a_{i'}Xa_{j'}$, $i\neq i'$, $j\neq j'$. Then $\phi a_j=Zy$, $\phi a_{j'}=Zy'$ for a word $Z$ and two letters $y\neq y'$, and by definition of the orders $y<_D y'$ if and only if $a_j<_{Dw} a_{j'}$. One of $(\phi a_i, \phi a_{i'})$ might be a strict suffix of the other, but each occurrence of $a_iXa_j$ or $a_{i'}Xa_{j'}$ in $L_v$ is preceded by an arbitrary long word, thus at least one occurrence of $\phi (a_iXa_j)$ in $L_{\phi v}$ is preceded by a word $U_{j_1}...U_{j_h}$ long enough
 to end with $w$, thus is preceded by $w$, and the same is true for at least one occurrence of $\phi (a_{i'}Xa_{j'})$ in $L_{\phi v}$. Thus $w\phi (a_iXa_j)$ and $w\phi (a_{i'}Xa_{j'})$ are in $L_{\phi v}$, and, as $i\neq i'$, there exist suffixes $xY$ of $w\phi a_i$ and $x'Y$ of $w\phi a_{i'}$ for letters $x\neq x'$, and by definition of the orders $x <_A x'$ if and only if $a_i<_{Aw} a_{i'}$. Thus we have proved $Y\phi (X)Z$ is bispecial in $L_{\phi v}$, and  thus by the hypothesis satisfies the order condition for $(<_{D},<_{A})$, hence $X$ satisfies the order condition for $(<_{Dw},<_{Aw})$.\\\

Thus we conclude by Theorem \ref{t1}.  A similar reasoning applies to prefix return words, considering words  of the form  $U_{i_1}\cdots U_{i_r}w$  instead of  $wU_{i_1}\cdots U_{i_r}$. \end{proof}
The proof of Theorem \ref{t2} is indeed a rather straightforward application of the definition 
of the derived orders, but it uses strongly the fact that the $U_i$ are return words: in a subword of $U_{i_1}\cdots U_{i_r}$, we know where are all the occurrences of $w$, and if there are none  the considered subword is in  $L$. Note that, because of the result mentioned  after Definition \ref{sy}, the integer $s$ is at most $k$.

This theorem applies immediately to the natural codings of   a (generalised)  interval exchange transformation $T$, and in that case, under a mild extra condition, we can identify the pair of derived orders. 

\begin{corollary}\label{clie} Let $T$ be  a generalised interval exchange transformation on  $k$ intervals. Let $w$ be a word in $L(T)$, Let $U_1, \ldots ,U_s$ be the
suffix (resp. prefix) return words to $w$ in $L(T)$.  Let $\phi$ be the morphism from $\mathcal A'=\{a_1, \ldots ,a_s\}$ to $\mathcal A^{\star}$ defined by $\phi a_i=U_i$. 
Then  for any word $v$ on $\mathcal A'$,  $\phi v$ clusters for the pair of orders of $T$ on $\mathcal A$  if and only if $v$ clusters for the pair of derived orders with respect to $w$ on $\mathcal A'$. 

If $T_{[w]}$, the induced map of $T$ on $[w]$, exists, the pair of derived orders with respect to $w$ is  the pair of orders defining $T_{[w]}$.\end{corollary}

\begin{proof} If $(<_D,<_A)$ is  the pair of orders of $T$, we apply Theorem 19 of \cite{fhuz} and Theorem \ref{t2} above. Under the extra  hypotheses above, the definitions of the induced map and of the orders allows us to identify the derived orders. \end{proof}

Note that we can deduce the ‘‘if" direction of Theorem \ref{simple} above from Corollary \ref{clie}, as a single letter always clusters for any pair of orders.
 
 For images by morphisms of words of length $2$,  Example \ref{e3} below provides a counter-example. \\
 
 A modified version of Theorem \ref{t2} applies if we consider, for a word $v$, the non extendable language $F(v)$  instead of   $L_v$, and it will be used in Section 5 below.
 
 \begin{proposition}\label{t2r} With all the notations of Theorem \ref{t2} and suffix, resp. prefix, return words, for any word $v$ on $\mathcal A'$,  $F(w\phi v)$, resp. $F(\phi vw)$,  satisfies the order condition for   $(<_D,<_A)$ on $\mathcal A$  if and only if $F(v)$  satisfies the order condition for  $(<_ {Dw},<_ {Aw})$ on $\mathcal A'$.
\end{proposition}
\begin{proof} Same as Theorem \ref{t2}.\end{proof}

We look now at two questions about the tightness of the above results. The first one is to know whether the words we consider could cluster for more pairs of orders than the one we state. The answer is negative in general, as, given the language of an interval exchange transformation, under some hypotheses we can  determine all the possible order conditions it can satisfy.

\begin{proposition}\label{2ord}  Let $L$ be the natural coding of a minimal aperiodic interval exchange transformation, on $k$ intervals, with no connection and orders $(<_D,<_A)$, and $L'$ be a language on the same alphabet, satisfying the order condition for  $(<'_D,<'_A)$ There exists $N$ such that, if $L\cap L'$ contains a word of length at least $N$,  then either $<'_A$ is the same as $<_A$ and 
$<'_D$ is the same as $<_D$, or $<'_A$ is the reverse of  $<_A$, and $<'_D$ is the reverse of  $<_D$.
\end{proposition} 
\begin{proof} $L$ is uniformly recurrent, and it is known, see \cite{fhuz} for example, that $L$ is has no strong bispecial, and that every long enough word in $L$ has at most  two left extensions and two right extensions. Hence there exist  $k-1$ words $A_1$, \ldots , $A_{k-1}$, $A_i$ having two left extensions, $a_{i,1}A_i$ and $a_{i,2}A_i$, such that every long enough  left special word in $L$ has as a prefix one of the $A_i$.  Similarly, there exist $k-1$ words $B_1$, \ldots , $B_{k-1}$, $B_i$ having two right extensions, $B_ib_{i,1}$ and $B_ib_{i,2}$, such that every long enough  right special word in $L$ has as a suffix one of the $B_j$. Because $L$ satisfies the order condition and  has no connections, there are one $a_{c,e}$ and one $a_{d,f}$, corresponding to the extremities of the order $<_A$, which are $a_{i,x}$ for only one $i$, and one $b_{c',e'}$ and one $b_{d',f'}$, corresponding to the extremities of the order $<_D$, which are $b_{i,x}$ for only one $i$. Also, by aperiodicity, there is at least one bispecial word  $u$ in $L$ having some $A_l$ as a prefix and some $B_m$ as a suffix, and it has three bilateral extensions $aub$, $a=a_{l,1}$ or $a=a_{l,2}$, $b=b_{m,1}$ or $b=b_{m,2}$.

 We choose $N$ such that each word $w$ in $L$ of length at least $N$ contains each $a_{i,1}A_i$, $a_{i,2}A_i$, $B_ib_{i,1}$ and $B_ib_{i,2}$, and  all the bilateral extensions of $u$.
 
Suppose $w$ is also in $L'$. Then the $A_i$ are left special in $L'$, the $B_i$ are right special in $L'$, $u$ is bispecial in $L'$, and all the $a_{i,1}A_i$, $a_{i,2}A_i$, $B_ib_{i,1}$ and $B_ib_{i,2}$, and all the $aub$ in $L$, are also in $L'$. If $L'$ contains $xA_i$, the fact that there are   at least $k$ different  $a_{i,l}$, and the order condition for $L'$ imply that $x$ is one of the $A_{i,l}$; similarly, if $B_iy$ is in $L$, $y$ is one of the $b_{i,l}$.
 Also, for each $i$, $a_{i,1}$ and $a_{i,2}$ must be consecutive in  the order $<'_A$, $b_{i,1}$ and $b_{i,2}$ must be consecutive in the order $<'_D$ Thus there are only two possibilities for the order $<'_A$, one starting from $a_{c,e}$, progressing from there through letters used twice $a_{i,2}= a_{j,1}$ and neighbours  $a_{i,1}$ and $a_{i,2}$, and ending at 
$a_{d,f}$, and the reverse order. Thus $<'_A$ is the same as $<_A$ or its reverse, and similarly $<'_D$ is the same as $<_D$ or its reverse. But the extensions of $u$ are either $aub$, $a'ub'$,  and $aub'$ or $a'ub$, when $a<'_A a'$ and $b<'_D b'$ or $a>'_A a'$ and $b>'_D b'$, or else $aub'$, $a'ub$, and $aub$ or $a'ub'$, when $a<'_A a'$ and $b'<'_D b$ or $a>'_A a'$ and $b'>'_D b$. Thus  among the possible couples made with $<_D$ or its reverse and $<_A$ or its reverse, the couple $(<'_D <'_A)$ takes exactly  two values, one being deduced from the other by reversing both orders. As the resolution of $u$ allows the value $(<_D,<_A)$, 
this proves our proposition. \end{proof}

\begin{example}\label{e3} By \cite{fhz3}, we can build a standard 3-interval exchange transformation, with the pair of orders $(1<_D 2<_D 3,3<_A 2<_A 1)$ where $w=12$ has the three suffix return words $\phi a=12$, $\phi b=1312$, $\phi c=212$. Then the pair of orders of $T_{[w]}$ is
$(a<_{Dw} b <_{Dw} c, b <_{Aw} c<_{Aw} a) $, and $bc$ does not cluster for $(<_{Dw}, <_{Aw})$, nor does $\phi (bc)=13112212$ for $(<_D, <_A)$. Also, $\phi(bc)$ contains all the words of length $2$ of $L(T)$, namely $11$, $13$, $31$, $22$, $12$, $21$. these are enough to identify the orders, so $\phi (bc)$ cannot cluster for any pair of  orders other than  $(1<_D 2<_D 3, 3<_A 2<_A 1)$ or the pair of reversed orders, thus cannot cluster for any pair of orders. \end{example}

Another question is what happens for interval exchange transformations {\it with flips},  where in Definition \ref{dit} we allow $T$ on $I_e$ to be (strictly) either increasing or decreasing, or, in the standard case, to have slope $\pm 1$, as  many results of \cite{fhuz} as still true for them. However, the present results do not extend to interval exchange transformations with flips, and indeed these provide examples of return words which do not cluster, because of the following proposition and the existence of (standard) {\it minimal} examples, the earliest one being due to A. Nogueira \cite{nog}, for which every long enough word in $L(T)$ contains all words of length $3$. 

\begin{proposition}\label{fi} Let $T$ be a generalised interval exchange transformation  with flips  on $k\geq 2$ intervals indexed in  $\mathcal A$. Suppose $\beta_1<\gamma_1<\beta_2 <\ldots <\beta_{k-1}<\gamma_{k-1}$, and $T$ has slope $-1$ on some $[a]$, for a letter $a$ which is neither maximal for $<_D$ nor minimal for $<_A$. Then every word in $L(T)$ which contains all words of length $2$ and all the extensions $xay$, $x\in\mathcal A$, $y\in\mathcal A$, in  $L(T)$, does not cluster for any pair of orders. \end{proposition}
\begin{proof} If $v$ clusters, $L_v$ satisfies the order condition for some pair of orders $(<'_D, <'_A)$. Let $(<_D, <_A)$ be the orders defined by $T$, $\bar <_D$ and $\bar <_A$ the respective reversed orders. $L_v$  contains all words of length $2$ of $L(T)$, which are $a_1b_1$ and  $a_jb_{j-1}$, $a_jb_j$, $2\leq j\leq k$, with $a_1<_A .... a_k$ and  
$b_1<_D .... b_k$ .As in the proof of Proposition \ref{2ord}, looking at  the bispecial empty word, we get that hat $(<'_D, <'_A)$ must be either ($<_D,<_A)$ or $(\bar <_D, \bar <_A)$. 
But by Proposition 8 of \cite{fhuz}, the bispecial word $a$ satisfies the order condition for $(\bar <_D,<_A)$; it has three extensions and they are all in  $L_v$, thus, again as in the proof of Proposition \ref{2ord}, we get that $(<'_D, <'_A)$ must be either $(\bar <_D, <_A)$ or $(<_D, \bar <_A)$, which gives a contradiction. \qed \\

\section{Rich languages and conservation of orders}
In Theorem 2, the morphism $\phi$ sends one  clustering to another one.  These two clusterings might be the same  when $\mathcal A=\mathcal A'$.  Then it   is easy to modify $\phi$, by renumbering the $U_i$, to make it preserve the lexicographical (departure) order,  so that the orders $<_{Dw}$ and $<_D$ are the same. But, having done that, when will we get that  $<_{Aw}$ and $<_A$ are the same? This will depend on $w$, in a way that does not seem to have been much studied. In general, it is known that the permutation giving $<_{Aw}$ from $<_{Dw}$ is in the {\it Rauzy class} of the permutation giving $<_A$ from $<_D$, and that all permutations in that class can be reached \cite{rau}. 

We shall give a more precise result in the (always easier to handle) symmetric case.
For this, we shall use the  theory of {\it rich languages}, defined and named in \cite{gjwz}, though this notion was also studied under different names in \cite{bdgz} and some similar ideas  appear already in \cite{bhnr}. The definition we give here is derived from 
Theorem~2.14 in \cite{gjwz}.

\begin{definition} A language $L$  is rich if  all complete return words to palindromes are palindromes.
\end{definition}

Throughout this section, we restrict ourselves to languages satisfying the symmetric order condition.   By Theorem~\ref{t1}, a word $w\in \mathcal A^+$ satisfies the symmetric order condition if and only if $w$ is a perfectly clustering word. For instance, taking $\mathcal A=\{a,b,c\}$ and $a<_D b<_D c$, $c<_A b<_A a,$ it is easily checked that the word $w=abac$ satisfies the symmetric order condition and is a perfectly clustering word. We now show that if a language $L$ satisfies the symmetric order condition, then $L$ is rich.   

\begin{lemma}\label{aric} Let $L$ be a language over an alphabet $\mathcal A$, satisfying the symmetric order condition for a pair of orders $(<_D, <_A),$ and let $w\in L$ be a complete first return word to a letter $a\in \mathcal A.$ Then $w$ is a palindrome.    
\end{lemma}

{\bf Proof} Assume $L$  satisfies the symmetric order condition for $(<_D, <_A).$ Write $w=aza$ with $z\in \mathcal A^*.$ Thus $a$ does not occur in $z.$ Set $u=az.$ Then $u$ is a prefix  return word to $a$ and thus by Theorem \ref{simple} $u$ is a perfectly clustering word. We now claim that the language $L_{\bar u}$  generated by $\bar u$ also satisfies the symmetric order condition for the same pair of orders $(<_D,<_A).$ In fact, suppose $axb, cxd \in L_{\bar u}$ with $a\neq c$ and $b\neq d.$ Then $b\bar x a, d\bar x c \in L_u$ and hence $$a<_A c \Leftrightarrow c<_D a  \Leftrightarrow d<_A b  \Leftrightarrow b<_D d.$$ By Theorem~\ref{t1}, the reverse word $\bar u$ is also a perfectly clustering word. However, since the Parikh vector of $u$ and $\bar u$ are the same, it follows that the Burrows-Wheeler transform of $u$ and $\bar u$ are the same. In other words $u=az$ and $\bar u =\bar z a$ are conjugate to one another. But as $a$ does not occur in $z,$ the only conjugate of $u$ ending in $a$ is $za.$ Thus $z=\bar z$ and  $w$ is a palindrome as required. \end{proof}

\begin{lemma}\label{bric} Let $a$ be a letter in $\mathcal A$, Let $L$ be a language on $\mathcal A$ satisfying the symmetric order condition for  $(<_D,<_A)$. Let  $U_1$, ... $U_s$ be 
all the suffix (resp. prefix) return words of $a$ in $L$, $1\leq s \leq k$. With the notations of Theorem \ref{t2}, 
if $w=a\phi z$, resp. $w=\phi z a$, and  $L_w$ satisfies the symmetric order condition for  $(<_D,<_A)$,  then $L_z$ satisfies the symmetric order condition for the derived orders with respect to $a$. \end{lemma}
\begin{proof}
We consider suffix return words; the proof for prefix return words is similar. We apply Proposition \ref{t2r}, but it remains to show  that $<_{Aa}$ is the reverse order of $<_{Da}$. By Lemma \ref{aric}, if $i\neq j$, then $aU_i$ and $aU_j$ are palindromes, 
thus $(U_i)_t=(aU_i)_u$ and $(U_i)_t=(aU_j)_u$ for $t$ and $u$ of Theorem \ref{t2}, thus we conclude as $<_{A}$ is the reverse order of $<_{D}$,.\end{proof}

\begin{theorem} \label{rich} Let $L$ be a language satisfying the  symmetric order condition. Then $L$ is rich. \end{theorem}

\begin{proof} Assume $L$ is a language satisfying the symmetric order condition for a pair of orders $(<_D, <_A).$ We prove by induction on $|w|$ that if $x\in L$ is a complete first return word to a non-empty palindrome $w,$ then $x$ is a palindrome. The base case of the induction, i.e., when $|w|=1,$ is precisely the content of Lemma~\ref{aric}.  So  let  $n\geq 2$ and  $x\in L$ be a complete first return word  to a palindrome $w$ of length $n.$ Let $a$ denote the first (and hence last) letter of $x.$ Then both $x$ and $w$ begin and end in $a.$ Thus we may write \[x=au_1au_2\cdots au_ma\] and \[w=av_1av_2a\cdots av_ka\] with $u_i,v_i \in L$ (possibly empty) not containing the letter $a$ and $k,m\geq 1.$ By Lemma~\ref{aric}, the $u_i$ and $v_i$  are each palindromes. Let $\mathcal A'=\{au_i\,|\, 1\leq i\leq m\}$ denote the set of all prefix  return words to $a$ occurring in $x.$  

Let $xa^{-1}$ (respectively $wa^{-1})$ denote the prefix of $x$ (respectively $w)$ obtained by deleting the last letter $a.$ Then $ xa^{-1}$ (respectively $wa^{-1})$ is a concatenation of prefix return words to $a$ and hence  may be viewed as a word $x'$ (respectively $w')$ of length $m$ (respectively $k)$ over the derived alphabet $\mathcal A'$, i.e. the alphabet whose letters are the prefix return words to $a$. Moreover, $x'$ is a complete first return to $w'.$ In fact, as $x$ begins and ends in $w$ it follows that $x'$ begins and ends in $w'.$ Moreover,  an internal occurrence of $w'$ in $x'$ would give rise to an internal occurrence of $w$ in $x.$ We claim that $w'$ is a palindrome over the alphabet $\mathcal A'.$ In fact, as $w$ is a palindrome,  it follows that $\bar v_i= v_{k+1-i}$ for 
each $1\leq i\leq k.$  Hence, $av_i=a\bar v_i=av_{k+1-i}$ for each $1\leq i\leq k$ and hence $w'$ is a palindrome. 
Also $n=|w|\geq k+1$ (with equality if each $v_i$ is empty) and hence $|w'|=k<n.$ Finally, by Lemma~\ref{bric},  $x'$ satisfies the symmetric order condition on the derived alphabet $\mathcal A'.$ Thus, by induction hypothesis, we have that $x'$ is a palindrome on the derived alphabet $\mathcal A',$ i.e., $au_i=au_{m+1-i}$ for each $i=1,\ldots m.$  Thus \[\bar x= a\bar u_ma \cdots a\bar u_2a\bar u_1a=au_ma\cdots au_2au_1a=au_1au_2a\cdots au_ma=x\] as required. 
\end{proof}

Note that  the converse to Theorem~\ref{rich} is false even in the binary case. For instance, the period doubling word, defined as the fixed point of the morphism $0\to  01$, $1 \to 00$ is uniformly recurrent and rich \cite{bal}.  However the period doubling word does not satisfy the symmetric order condition as it contain $000$ and $101$ as factors. \\

Special cases of Theorem \ref{rich} were already known. For example,  if $T$ is a standard symmetric interval exchange transformation without connections, $L(T)$ is a rich language. This comes from Theorem 1.1 of  \cite{bdgz},  Theorem 2.14 of \cite{gjwz}, and the computation of the {\it palindromic complexity} $PC(n)$ of $L(T)$ in  Theorem 4.1 of \cite{pel}. Namely, $PC(n)$ is the number of distinct palindromes of length $n$ in a language, and, for such $L(T)$, $PC(n)$ is $1$ for $n$ even, $k$ for $n$ odd; the proof in \cite{pel} re-discovers the {\it Weierstrass points} associated to $T$, see for example \cite{ll}, proving that all palindromes are in their trajectories. This leads us to study the palindromic complexity in the  more general case of Theorem \ref{rich}.

\begin{proposition}\label{palco} Let $L$ be a language on $k$ letters satisfying a symmetric order condition. Then its palindromic complexity is at most  $1$ for $n$ even, $k$ for $n$ odd.

Furthermore, if $L$ is generated by an infinite or bi-infinite word, is closed under reversal and has no connections, its palindromic complexity is $1$ for $n$ even, $k$ for $n$ odd.
\end{proposition}
\begin{proof} There are one palindrome of length $0$ and $k$ palindromes of length $1$. Then, when we go from $n$ to $n+2$, $awa$ and $bwb$ cannot both exist if $a\neq b$ as this contradicts the symmetric order condition, which gives the first assertion. 

By Theorem 1.1 of  \cite{bdgz}, for a rich  language generated by an infinite word and  closed under reversal, $PC(n)+PC(n+1)=p(n+1)-p(n)+2$, where  $p(n)$ is the usual factor complexity, counting the number of different factors of length $n$ in $L$; as $L$ has no connection $p(n+1)-p(n)=k-1$, see for example Corollary 4 of \cite{fhuz}, and this together with the first assertion gives the second one. \end{proof}

The second assertion of Proposition \ref{palco} applies to Examples 4 and 9 of \cite{fhuz}, which are outside the scope of \cite{pel}. But it does not always hold if $L$ has connections.

\begin{example}\label{e} Let $L$ be generated by the bi-infinite word $(1312)^{\omega}$. It satisfies the symmetric order condition for  $(1<_D 2<_D , 3<_A 2<_A 1)$, is closed under reversal, but its palindromic complexity is $0$ for even  $n>0$, $2$ for odd $n>1$.  Indeed, as will result from Section 5 below, $L$ is
 the natural coding of a standard  interval exchange transformation, but with connections.  \end{example}
 
 We look now at a language which is not closed under reversal, as it contains the word $3122$ but not $2213$.
 
\begin{example}\label{e2} Let $L$ be generated by the bi-infinite words  $(1312)^{\omega-}(212)^{\omega+}$ and $(3121)^{\omega-}(221)^{\omega+}$. It  satisfies the symmetric order condition for  $(1<_D 2<_D 3, 3<_A 2<_A 1)$, has no connection, and  its palindromic complexity is indeed $1$ for 
  $n$ even, $3$ for $n$ odd. We remark that all its palindromes come from the sublanguage $L'$ 
 generated by the bi-infinite words  $(1312)^{\omega}$ and $(221)^{\omega}$, which is  the natural coding of a standard  interval exchange transformation with connections, but still has maximal palindromic complexity.
  \end{example}
 
In Example \ref{e2}, the connections in $L'$ are only on the bispecials $12$ and $21$, which are reverses of each other, and thus their effect is to forbid only words of the form $x12y$ or $x21y$, for letters $x$, $y$, thus the connections do not reduce the number of palindromes. The theory in  \cite{fz1}, see
the proof of Proposition \ref{pfls} below, suggests that this phenomenon does generalise, thus we make the following conjecture.

\begin{conj}  A language on $k$ letters satisfying a symmetric order condition and with no connection has palindromic complexity $1$ for $n$ even, $k$ for $n$ odd.
\end{conj}

We proceed now to determine in which cases the derived orders coincide with the original ones, which is the case (up to a renumbering of the return words) if and only if   $<_{Aw}$ is the reverse of $<_{Dw}$,
. We first state the result in its fullest generality; let the $U_i$, $<_D$, $<_A$, $<_{Dw}$, $<_{Aw}$ be as in Theorem \ref{t1}.

\begin{proposition}\label{sym} 
Let $L$ be a language satisfying the symmetric order condition for $(<_D, <_A)$,  if $w$ in $L$ is a palindrome or if the shortest bispecial word containing $w$ exists and is a palindrome, then  $<_{Aw}$ is the reverse of $<_{Dw}$,\\
Let $L$ be a recurrent  language satisfying the symmetric order condition for  $<_D, <_A)$,  if the shortest bispecial word containing $w$ exists and is not a palindrome, then  $<_{Aw}$ is not the reverse of $<_{Dw}$.
 \end{proposition}
\begin{proof}
First we notice that if $w$ is extended uniquely  right and left to get a bispecial word $w'$, the derived orders with respect to $w$ and $w'$ are the same.

Under the hypothesis of the first assertion, if $w$ is a palindrome or $w'$ is a palindrome, by Theorem \ref{rich} all the full return words of $w$, resp. $w'$, are palindromes.  Thus by the proof of Lemma \ref{bric} we conclude that $<_{Aw}$ is the reverse order of $<_{Dw}$, resp. $<_{Aw'}$ is the reverse order of $<_{Dw'}$.\\

Under the hypothesis of the second assertion,  $L(T)$ is a rich language by Theorem \ref{rich}, recurrent by Theorem 13 of \cite{fhuz}, and closed under reversal by Lemma 6 of \cite{fz3}. Thus $\bar w'$ is bispecial.  If $w'$ is not a palindrome, by Propositions 2.3 and 2.4 of \cite{bdgz}, every complete return word $v$ of $w'$ 
contains one occurrence of $\bar w'$, and, if $v=z\bar w'z'$, $z\bar w'$ and $\bar wz'$ are palindromes.

 Let $U_i$ be the suffix return words of $w'$ (a similar proof can be made by using the prefix return words); we use the notations of Theorem \ref{t1}. 
We choose two extensions $x\bar w'y$ and $x'\bar w'y'$ with $x\neq x'$ and $y\neq y'$. As $\bar w'$ satisfies the order condition for $(<_D,<_A)$, we suppose for example $x<_A x'$, $y<_A y'$. As $L$ is recurrent, we can extend $x\bar w'y$ to
the left until we get a word $X\bar w'y$ with $w'$ as a prefix, with the shortest possible $X$, then  extend $X\bar w'y$ to
the right until we get a word $X\bar w'Y$ with $w'$ as a suffix, with the shortest possible $Y$; thus $X\bar w'Y$ is a complete return word of $w'$. Thus it is a $w'U_i$, and  also it has exactly one occurrence of $\bar w'$. Similarly, we extend $x'\bar w'y'$ to $X'\bar wY'=w'U_j$. As $X\bar w'$, $\bar w'Y$, $X'\bar w'$, $\bar w'Y'$ 
are all palindromes, we get that, for the $t$ and $u$ in Theorem \ref{t2},  $(U_i)_t=x$, $(U_j)_t=x'$, $(w'U_i)_u=y$ and $(w'U_j)_u=y'$. As $x<_A x'$, we have $x' <_D x$, thus $a_j <_{Dw}a_j$; as $y<_D y'$, we have $y' <_A y$, thus $a_j <_{Aw}a_j$, thus $<_{wA}$ is
 not the reverse order of $<_{wD}$.\end{proof}
 
 Proposition \ref{sym} can be readily  translated into the vocabulary of   interval exchange transformations, the most interesting case being the following.

\begin{corollary}\label{cosym}
Let $w$ be a word in $L(T)$, for $T$ a  symmetric standard interval exchange transformation with no connection. Then $T_{[w]}$ is symmetric if  and only if the shortest bispecial word containing $w$ is a palindrome.
 \end{corollary}
\begin{proof}
In this case, $L(T)$ is aperiodic and (uniformly) recurrent, thus $v$ can always be extended to a bispecial, and we apply Proposition \ref{sym}. \end{proof}

\section{Non-clustering words and order condition}
 An  application of the notion of clustering is:
 
 \begin{proposition} For any primitive word $w$,  the following properties are equivalent 
 \begin{itemize}
  \item[(i)] $ww$ is produced by
a standard  interval exchange transformation with orders   $(<_D,<_A)$,
   \item[(ii)] $ww$ is produced by 
a generalised  interval exchange transformation with orders   $(<_D,<_A)$,
 \item[(iii)] $w$ clusters for the pair of  orders   $(<_D,<_A)$.
 \end{itemize}
\end{proposition}
\begin{proof} (i) implies (ii) trivially. (ii) implies (iii) by Proposition 8 of \cite{fhuz}, as it  implies  $ww$ is in some $L(T)$ which satisfies the order condition, Theorem \ref{t1} and Lemma \ref{useful}. (iii) implies (i) by Theorem 1 of \cite{fz4}.\end{proof}

Thus we can wonder if $(i)$ is still equivalent to $(ii)$  for non-clustering words, when we replace $vv$  by just $v$. Thus we must look now at properties of finite words satisfying some kind of order conditions, in the same way  as \cite{fhuz}  looks at properties of whole languages satisfying order conditions.
 
\begin{proposition}\label{pfl}
Let $W$ fe a finite set of words. The non extendable language $F(W)$ satisfies the order condition for  $(<_D,<_A)$ if and only if $W$ is included in an (extendable) language satisfying the order condition for $(<_D,<_A)$, with no connection.
\end{proposition}
\begin{proof}
The ``if" direction is immediate. Now, we start from  $F(W)$, and fix the two orders 
$<_D$ and $<_A$. We introduce the {\it ordered extension graph} of a bispecial word $w$ in a language $L$ equipped with a pair of orders, which are a small modification of the extension graphs of \cite{b+} and \cite{d} on an alphabet with orders. We put on a line the letters $x$ such that $xw$ is in $L$, ordered from left to right according to $<_A$, and on a line above the letters $y$ such that $wy$ is in $L$, ordered from left to right according to $<_D$, and draw  an edge from $x$ (on the line below)  to $y$ (on the line above) if $xwy$ is in $L$; then $w$ satisfies the order condition if and only if no two edges in its ordered extension graph intersect except possibly at their endpoints. 

Given the ordered extension graph of $w$ in $L$, we define its  {\em saturate}:  if it is possible to add an edge without intersecting the existing ones, we add the one with the leftmost possible origin, and iterate the process with 
the new graph, and iterate again until we cannot add any extra edge. The resulting graph is the extension graph of $w$ in another (larger) language, and we say this extension graph is saturated; in it $w$ still satisfies the order condition for $(<_D,<_A)$. Note that in the saturated graph there is at least one edge starting from every $x$ on the bottom line: otherwise, we look at every edge starting left of $x$, call $y_0$ the rightmost of their arrival points, 
and we can draw an extra edge from $x$ to $y_0$. Also, there is at least one edge arriving at every $y$ on the top line: otherwise, we look at every edge arriving left of $y$, call $x_0$ the rightmost of their departure points, 
and we can draw an extra edge from $x_0$ to $y$. And, if the ordered extension graph of $w$ is saturated, then $w$ has no  connection: otherwise, with $a,a',b,b'$ as in Definition \ref{con}, we can add an extra edge from $a$ to $b'$ or from $a'$ to $b$.\\

Let $\mathcal A$, of cardinality $k$, be the set of all letters in  $F(W).$ We  take the saturate of the ordered extension graph of the empty bispecial in $F(W)$ and define $ W_2$ to be the set of words $xy$ such that there is an edge form $x$ to $y$ in the saturated graph; 
$W_2$ contains the words of length $2$ of $F(W)$. Suppose now we have defined $W_n$ and it contains the words of length $n$ of $F(W)$, we define now $W_{n+1}$: for a word $w$ of length $n-1$ in $W_n$, we look at  all its extensions $xwy$, $xw\in W_n$, $wy \in W_n$; 
if $w$  is not bispecial, we put in $W_n$ all its extensions $xwy$; if it is bispecial, we draw
 its ordered extension graph in $F(W)$, defined in the usual way  if $w$ is in $F(W),$ even if it is not  bispecial in $F(W)$, and without any edge if $w$ is not in $F(W)$, take its saturate, and put in $W_{n+1}$ all the $xwy$ such that there is an edge from $x$ to $y$ in the saturated graph. Note that for $n$ large enough there is no word of length $n$ in  $F(W)$, thus the extension graphs we have to saturate have no edges; taking their saturate gives  just one of the possible way to resolve the bispecials  without contradicting the order condition.\\

Then the language $\cup_{n\geq 2} W_n$ is extendable and has no connection by the above considerations on saturated ordered extension graphs, satisfies the order condition and contains $F(W)$.\end{proof}

\begin{example}\label{e4}  $W=\{11,22\}$; $F(W)$ satisfies the order condition for  $(1<_D 2,1<_A 2)$. The saturation of the ordered extension graph of the empty word $\varepsilon$ gives
 $W_2=\{11,12,22\}$, and then there is no bispecial word of length at least $1$, $W_2$ determines the language  $\cup_{n\geq 2} W_n$, which is generated by the bi-infinite word $1^{\omega-}2^{\omega+}$, as in  Example 3 of \cite{fhuz}. The picture shows the process of saturation for the bispecial empty word.\end{example}

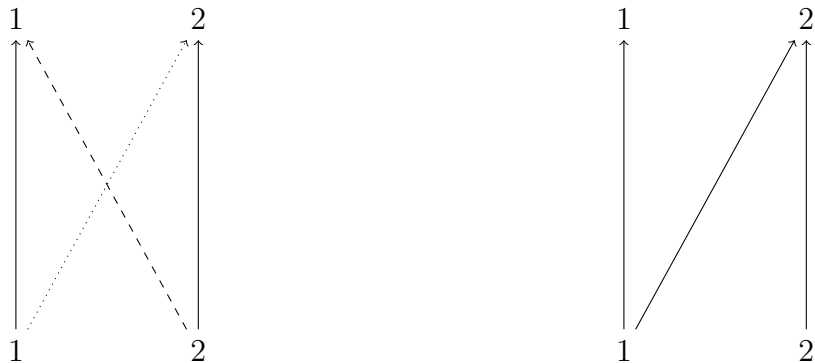
\begin{center}
\begin{tikzpicture}[scale = 8][every text node part/.style={align=center}]

\node (E) at (.5,-.05) {$1$};
\node (F) at (.8,-.05) {$2$};
\node (G) at (.5,.50) {$1$};
\node (H) at (.8,.50) {$2$};

	\draw[->] (E) edge node[auto] {} (G);		
	\draw[->] (F) edge node[auto] {} (H);		\draw[->,dashed] (F) edge node[auto] {} (G);	
	
			 \draw[->,dotted] (E) edge node[auto] {} (H);

\node (A) at (1.5,-.05) {$1$};
\node (B) at (1.8,-.05) {$2$};
\node (C) at (1.5,.50) {$1$};
\node (D) at (1.8,.50) {$2$};

	\draw[->] (A) edge node[auto] {} (C);		
	\draw[->] (B) edge node[auto] {} (D);		
	
		\draw[->](A) edge node[auto] {} (D);

\end{tikzpicture}
\end{center}

\begin{corollary} $F(W)$ satisfies the order condition for  $(<_D,<_A)$ if and only if $W$ is produced by
a generalised interval exchange transformation with the same pair of  orders.\\
With the notations of Theorem \ref{t2} for the morphism $\phi$, if the $U_i$ are suffix (resp. prefix) return words,  $w\phi (v)$ (resp. $\phi (v)w$)  is produced by a generalised  interval exchange transformation with the same orders as $T$ if and only if $v$ is produced by a generalised interval exchange transformation with the pair of orders $<_{Dw},<_{Aw}$.\end{corollary}
\begin{proof} The first assertion is deduced from Proposition \ref{pfl} through Theorem 19 of  \cite{fhuz}, the second one comes from Proposition \ref{t2r}.\end{proof}

 We want now to improve the nontrivial direction of Proposition \ref{pfl}  to get something better than a natural coding of a generalised interval exchange transformation. The following  result could also be deduced from Proposition \ref{pfl} 
and the proof of Theorem 19 of  \cite{fhuz}, but the present proof is self-contained and uses only words of bounded length. Its proof uses the fact that resolving bispecials in the language of an  interval exchange transformation is equivalent to ordering the images of the discontinuities.

\begin{proposition}\label{pfl2}
The non extendable language $F(W)$ satisfies the order condition for  $(<_D,<_A)$ if and only if $W$ is included in a grouped coding of
an affine interval exchange transformation with the same pair of  orders.\end{proposition}
\begin{proof}
The ``if" direction is immediate. 

In the other direction, let $\mathcal A$, of cardinality $k$, be the set of all letters in  $F(W).$ Starting from $F(W)$, we build the $W_n$ as in Proposition \ref{pfl}, but only   when $n\leq N$, $N$ being  the maximal length of the words in $W$.  We have  $\mathcal A=\{b_1<_D\ldots <_D b_k\}$,  $\mathcal A=\{a_1<_A\ldots <_A a_k\}$. \\
 Note that if $w$ is in $W_{n-2}$  and we look at its extensions in $W_n$, it can be followed by $b_j,,\ldots ,b_{j+m}$ for some $j$ and $m,$, preceded by $a_i, \ldots ,a_{i+h}$ for some $i$ and $h$ (if there is a gap, this contradicts the order condition and saturation at a  previous stage). \\

We recall  some properties of words of length $n$ in the language $L(T)$ for an interval exchange transformation $T$ with no connection. Such a word $w$ corresponds to  the cylinder $[w]$ for $T$, which is an interval $J_w=I_{w_1}\cup T^{-1}I_{w_2}.\cup \ldots \cup T^{-n+1}I_{w_n}$ whose endpoints are among the $T^p\gamma_i$, $0\leq p\leq n-1$.   
A word $w$ of  length $n$ is left 
special whenever $J_w$ contains a $\beta_l$, right special whenever $J_w$ contains a $T^{-n}\gamma_l$. Suppose  $w$ is bispecial; it can be preceded by $a_i<_A\ldots <_A a_{i+h}$, and followed by $b_j<_D\ldots <_D b_{j+m}$. Then $J_w$ contains $T^{-n}\gamma_l$, $j \leq l\leq j+m-1$
and $\beta_l$, $i \leq l\leq i+h-1$. By the geometry of $T$, $a_ib_j$ is in $L(T)$; then $a_ib_{j+1}$ is in $L(T)$ if $T^{-n}\gamma_j<\beta_i$, $a_{i+1}b_j$ is in $L(T)$ if $T^{-n}\gamma_j>\beta_i$. Continuing from left to right, we get from the resolution of $w$ a unique (strict) ordering
 of the concerned $T^{-n}\gamma_l$ and $\beta_l$ (for the usual order, left to right, on the interval); we call it the $(L(T), w)$-ordering. It can be described by integer intervals $Z_l$, $i-1\leq l\leq i+h$, some of them possibly empty, with $l_1<l_2$ when $l_1\in Z_l$, $l_2\in Z_{l+1}$,
 whose disjoint union is the integer interval $[j,j+m-1]$, such  that $\beta_l<T^{-n}\gamma_{l'}< \beta_{l+1}$ when $l'$ is in $Z_l$, $i \leq l\leq i+h-1$, $T^{-n}\gamma_{l'}< \beta_{i}$ when $l'$ is in $Z_{i-1}$, $T^{-n}\gamma_{l'}> \beta_{i+h}$ when $l'$ is in $Z_{i+h}$.

 Let $w$ be a bispecial word in $W_{n-2}$, and we look at its extensions  in $W_n$; then  it can be followed by $b_j,\ldots ,b_{j+m}$ for some $j,m$, preceded by $a_i, \ldots ,a_{i+h}$ for some $i,h$ (if there is a gap, this contradicts the order condition and saturation at a  previous stage). If $L$ was an $L(T)$ for an interval exchange transformation $T$ with no connection, $W_n$ would define an $(L(T),w)$ ordering described by  integer intervals $Z_l$, $i-1\leq l\leq i+h$, as above.
 This ordering depends only on $W_n$ because of the saturation of the extension graph (there has to be an edge form $a_i$ to $b_j$, then either an edge from $a_i$ to $b_{j+1}$ or an edge from $a_{i+1}$ to
$b_j$, 
and so  on), we call it the $(W_n,w)$ ordering.  \\

 We shall now build interval exchange transformation $T_n$ so that the set of words of  length $n$ of  $L(T_n)$ is $W_n$. We consider them  as generalised interval exchange transformations with the same defining intervals: indeed, the $T_n$  are piecewise affine on their (common) defining intervals, and are also affine interval exchange transformation but with a growing number  of defining intervals.  We shall denote by $J_{n,w}$ the interval corresponding to the cylinder $[w]$ for $T_n$. At the beginning we choose
points $0=\gamma_{0}<.\ldots < \gamma_{k}=1$, define  $I_{b_l}=[\gamma_l,\gamma_{l+1})$ for $0\leq  l \leq k-1$.

Then we  choose  points $0=\beta_0<\ldots < \beta_{k}=1$, such that the natural ordering (from left to right) of the points  $\gamma_l$, $1\leq l\leq k-1$
and $\beta_l$, $1\leq l\leq k-1$ is the  same as the $(W_2, \varepsilon)$ ordering. We  define $T_2$ as the affine  interval exchange transformation for which  $I_{b_l}=[\gamma_l,\gamma_{l+1})$ for $0\leq i\leq k-1$, and $T_2I_{a_l}=[\beta_l,\beta_{l+1})$ for $0\leq i\leq k-1$. Then the set of words of length $2$ of $L(T_2)$ is $W_2$. 

To build $T_3$, we look at words of length $1$ and their bilateral extensions in $W_3$. If $w$ is not bispecial, we put $T_3=T_2$ on $J_{2,w}$. If $w$ is bispecial, preceded by $a_i<_A\ldots <_A a_{i+h}$, and followed by $b_j<_D\ldots <_D b_{j+m}$, its extensions in $L(T_2)$ are determined by the $(L(T_2),w)$ ordering
 while its extensions in $W_3$ are determined by the $(W_3, w)$ ordering, given by integer intervals $Z_l$ as above, and these may give different extensions. So we shall modify $T_2$ into $T_3$ such that  the $(W_3, w)$ ordering and the $(L(T_3), w)$ ordering  give the same extensions; for this, we choose arbitrarily some new points
  $\gamma_{i,3}$ in $J_{2,w}$, 
 such that 
 \begin{itemize}
 \item $\gamma_{l,3}<\gamma_{l',3}$ whenever $T_2^{-1}\gamma_l<T_2^{-1}\gamma_{l'}$,
 \item   $\beta_l<\gamma_{l',3}< \beta_{l+1}$ when $l'$ is in $Z_l$, $i\leq l\leq i+h-1$, $\gamma_{l',3}< \beta_i$ when $l'$ is in $Z_{i-1}$, $\gamma_{l',3}> \beta_{i+h}$ when $l'$ is in $Z_{i+h}$.
\end{itemize}
 When $l$ is such that $T_2^{-1}\gamma_l$ is in no $J_{2,w}$ for any bispecial $w$, we put $\gamma_{l,3}=T_2^{-1}\gamma_l$. Thus we  have added $k-1$ points $\gamma_{i,3}$. And on $J_{2,w}$ we replace the affine $T_2$ by  a piecewise affine $T_3$  such that  $T_3J_{2,w}=T_2J_{2,w}$ 
  and $T_3\gamma_{i,3}=\gamma_l$; thus  the set of words of length $3$ of $L(T_3)$ is $W_3$. \\

Similarly, suppose we have built $T_{n-1}$. To build $T_n$, we look at words of length $n-2$, in $W_{n-2}$  or equivalently in  $L(T_{n-2})$ or  $L(T_{n-1})$,  and their bilateral extensions in $W_n$. If $w$ is not bispecial, we put $T_n=T_{n-1}$ on $J_{n-1,w}$. If $w$ is bispecial, preceded by $a_i<_A\ldots <_A <a_{i+h}$, and followed by $b_j<_D\ldots <_D b_{j+m}$, we modify $T_{n-1}$ into $T_n$ such that  the $(W_n, w)$ ordering, given by the integer intervals $Z_l$, and the $(L(T_n), w)$ ordering  give the same extensions; for this, we create new points $\gamma_{i,n}$ in $J_{n-1,w}$, such 
that 
 \begin{itemize}
 \item $\gamma_{l,n}<\gamma_{l',n}$ whenever $T_{n-1}^{-1}\gamma_{l,n-1}<T_{n-1}^{-1}\gamma_{l',n-1}$ 
  \item   $\beta_l<\gamma_{l',n}< \beta_{l+1}$ when $l'$ is in $Z_l$, $i\leq l\leq i+h-1$, $\gamma_{l',n}< \beta_i$ when $l'$ is in $Z_{i-1}$, $\gamma_{l',n}> \beta_{i+h}$ when $l'$ is in $Z_{i+h}$.
\end{itemize}
When $l$ is such that $T_{n-1}^{-1}\gamma_{l,n-1}$ is in no $J_{n-1,w}$ for any bispecial $w$, we  put $\gamma_{l,n}=T_{n-1}^{-1}\gamma_{l,n-1}$. And on $J_{n-1,w}$ we replace $T_{n-1}$, which is affine on this set, by  a piecewise affine $T_n$  such that 
 $T_nJ_{n-1,w}=T_{n-1}J_{n-1,w}$  and $T_n\gamma_{l,n}=\gamma_{l,n-1}$. Thus we  have added $k-1$ points $\gamma_{l,n}$, and, by the induction hypothesis  that $T_{n-1}^{n-3}\gamma_{l,n-1}=\gamma_l$, we have   $T_n^{-n+2}\gamma_l=\gamma_{l,n}$, which ensures that  the set of words of length $n$ of $L(T_n)$ is $W_n$.

Thus the proposition is proved with the interval exchange transformation $T_N$. \end{proof}

 Note that the final $T_N$ may have connections, though this will not happen on bispecials of length at most $N-2$.  Proposition \ref{pfl2} might possibly be improved by omitting the word ``grouped"l, as in  all the examples  in \cite{fhuz} for which the language is not the natural coding of an affine interval exchange transformations (Examples 4 and 9 of that paper) 
 every finite set of words is included in the natural coding of an affine  interval exchange transformation, and this is also the case in Example \ref{e5} below.
 
 \begin{conj}   $F(W)$ satisfies the order condition for  $(<_D,<_A)$ if and only if $W$ is in the natural coding of
an affine interval exchange transformation with the same orders.\end{conj}

But what about getting a natural coding of a standard interval exchange transformation? 
 
\begin{question} Let $F(W)$ satisfies the order condition for  $(<_D,<_A)$. Under which condition is $W$   produced by
a standard interval exchange transformation with the same orders? \end{question}

In the non-symmetric cases, we do not have any explicit condition ensuring that $W$ is in the natural coding of a standard interval exchange transformation, but we can follow the algorithm of the proof  of Theorem 13 of \cite{fhuz}: we define parameters which are the tentative measures of the intervals corresponding to the words of length $1$, and build a language by resolving the bispecials in a way satisfying the order condition;  we use the initial parameters to compute the tentative measures of the intervals corresponding to the words of length $n$, as in the proof mentioned above. If there exists a set  of parameters for which these measures are positive for all $n\leq N$, the maximal length of a word in $W$, then $W$ is in the natural coding  of a standard interval exchange transformation. Otherwise, there is an obstruction for some  $n\leq N$.  

\begin{example} $w=3322$, $F({w})$  satisfies the order condition for  $(1<_D 2<_D 3,2<_A 3<_A 1)$, and we can find numbers $(L_1,L_2, L_3)$ such that the natural coding of a $3$-interval exchange transformation with these  orders such that $[i]$ has length $L_i$, $i=1,2,3$, contains $w$. \end{example}

\begin{example} $w=3311$. $F({w})$  satisfies the order condition for  $(1<_D 2<_D 3,2<_A 1<_A 3)$  (which corresponds to a non-primitive permutation, in the sense of Section 2). If we want to produce $w$ by a standard $3$-interval exchange transformation with these orders, there is an immediate obstruction: the interval $[31]$ has length zero, thus $w$ is not produced by any standard interval exchange transformation with these orders. But, using the proof of Proposition \ref{pfl}, we can see that $w$ is included in the language  generated by the bi-infinite words $3^{\omega-}1^{\omega+}$ and $3^{\omega-}21^{\omega+}$, which is a natural coding of a generalised (indeed affine) interval exchange transformation with these orders.  
Note that of course we could use the previous example by exchanging $1$ and $2$, but we can check that $3311$ is not produced by any standard $3$-interval exchange transformation with $1<_D 2<_D 3$, nor of any standard $2$-interval exchange transformation with $1<_D 3$. If we go to $4$ intervals, $w$ satisfies the order condition for  $(1<_D 2<_D 3 <_D 4,4<_A 2<_A 1<_A 3)$, corresponding to a primitive permutation, but  if we want to produce $w$ by a standard $4$-interval exchange transformation with these orders, we find another kind of  obstruction: the interval $[1]$ has to be strictly longer than $T[1]$. \end{example}

Unsurprisingly,  the symmetric case is the one in which we can prove general results. For them, we rely heavily on the machinery of  \cite{fz1}.
\begin{proposition}\label{pfls}
  $W$ is produced by 
a standard symmetric  interval exchange transformation with orders $(<_D,<_A)$ and no connection if and only if $F(W)\cup F(\bar W)$ satisfies the order condition for  
$(<_D,<_A)$.\end{proposition}
\begin{proof}
The ``only if" part comes from the fact, proved in  Lemma 6 of  \cite{fz3}, that 
the language of a symmetric interval exchange transformation is closed under reversal.

Let $\mathcal A$, of cardinality $k$, be the set of all letters in  $F(W).$ In \cite{fz1}, we described the language  of a standard symmetric interval exchange transformation through a sequence of objects called {\it (composite) trees of relations} and {\em instructions}, the instructions corresponding to resolution of bispecial words; then Theorems 2.9 and 3.1 of that paper give conditions on such sequences ensuring that they generate the language of a standard symmetric interval exchange transformation with no connection. We shall build such a sequence here.\\

We  take the  saturate of the ordered extension graph of the empty bispecial in $F(W)\cup F(\bar W)$ and define the non extendable language $Y_1$ to be the set of factors of words $xy$ such that there is an edge form $x$ to $y$ in the saturated graph; 
$Y_1$ contains the words of length $2$ of $F(W)\cup F(\bar W)$, and, because of the definition of the saturate, is closed under reversal. Then the bispecials $w_{1,i}$ of length $1$ of $Y_1$ can be described by a (composite in general) tree of relations $\Gamma_1$.

At stage $n$, suppose we have defined a non extendable language $Y_n$, closed under reversal, whose bispecial words are  those of $Y_{n-1}$ plus some $w_{i,n}$ (with possibly, for some $i$,  $w_{i,n}=w_{i,n-1} \in Y_{n-1}$),  which can be described by a composite tree of relations $\Gamma_n$, and we  define an instruction $\iota_n$ on $\Gamma_n$ in the following way: if $w_{i,n}$ is  in $F(W)\cup F(\bar W)$, we choose a resolution of each $w_{i,n}$ in our usual way, by using the saturate of the ordered extension graph of $w_{i,n}$ in $F(W)\cup F(\bar W)$, and define $\iota_n$  on $w_{i,n}$ so that it corresponds to this resolution; note that the fact that $Y_n$ is closed under reversal, the order condition for $F(W)\cup F(\bar W)$, and the definition of the saturate, imply that $\iota$ is the same on  $w_{i,n}$ and  $\bar w_{i,n}$. For the    $w_{i,n}$ which are not  in $F(W)\cup F(\bar W)$,  we choose any instruction on them which  is the same on  $w_{i,n}$ and  $\bar w_{i,n}$.

Then, as in \cite{fz1},  $\Gamma_n$ and $\iota_n$ define a new composite tree of relations $\Gamma_{n+1}$, a set of words $w_{i,n+1}$, with possibly $w_{i,n+1}=w_{i,n}$ for some $i$, and a new non extendable language $Y_{n+1}$ made with the factors of all the $w_{i,n+1}$, and we continue the recursion. There exists $M$ such that, when $w_{i,n} \neq w_{i,n+1}$ for  $M$ different integers $n_1<...  n_M$, then $w_{i,n_M}$ is too long to be  in $F(W)\cup F(\bar W)$. Thus we can ensure, as in Lemma 2.5 and Proposition 2.6 of \cite{fz1}, that our construction is never stuck, and that
we can  choose the infinite sequence of $\Gamma_n$ and $\iota_n$ so that it is admissible in the sense of Theorems 2.9 and 3.1 in \cite{fz1}. Also, for $n$ large enough, $Y_n$ contains $F(W)\cup F(\bar W)$.
Thus we have built a language of a standard symmetric interval exchange transformation with no connection, which contains $F(W)\cup F(\bar W)$. \end{proof}

\begin{example}\label{e5b} If $W=\{32, 12\}$, $F(W)$ satisfies the symmetric order condition for  $(1<_D 2<_D 3,3<_A 2<_A 1)$, but $F(W)\cup F(\bar W)$ does not satisfy it. $W$ is not included in the natural coding of any symmetric interval exchange transformation, but,  using the proof of Proposition \ref{pfl}, 
we can still put $W$ in
 the natural coding $L$ of a generalised, indeed affine, interval exchange transformation; this language is generated by the bi-infinite words  $(13)^{\omega-}2^{\omega+}$ and $(31)^{\omega-}2^{\omega+}$. Note that $L$   has no nonempty bispecial word, is 
 not recurrent, and not closed under reversal. $L$ is the simplest language we have found which satisfies the symmetric order condition with no connection but  contains a finite subset which  is not produced by any standard interval exchange transformation.  \end{example}

\begin{example}\label{e5} If $W=\{3122, 1212\}$, $F(W)$ satisfies the symmetric order condition for  $(1<_D 2<_D 3,3<_A 2<_A 1)$, but $F(W)\cup F(\bar W)$ does not satisfy it. $W$ is not included in the natural coding of any symmetric interval exchange transformation, but,  using the proof of Proposition \ref{pfl}, 
we can still put $W$ in
 the natural coding $L$ of a generalised, indeed affine, interval exchange transformation; this language is the one of Example \ref{e2},  generated by the bi-infinite words  $(1312)^{\omega-}(212)^{\omega+}$ and $(3121)^{\omega-}(221)^{\omega+}$. Note that $L$   has no bispecial word of length $>2$, is 
 not recurrent, and not closed under reversal.  \end{example}
 
 But when $W$ is made with a single word, we can prove more: the following theorem implies that, in the symmetric case, $w$ is produced by a
    generalised interval exchange transformation if and only if $w$ is produced by a standard interval exchange transformation.

 \begin{theorem} \label{c2} A word $w$ is produced by
a standard symmetric interval exchange transformation  with orders  $(<_D,<_A)$ and no connection  if and only if  $F({w})$ satisfies the order condition for the same pair of orders.\end{theorem}
\begin{proof}
We have only to prove the ``if"   direction. Suppose   $F({w})$ satisfies the symmetric order condition. If Proposition \ref{pfls}  can be applied, we have finished. Otherwise, $F(w)\cup F(\bar w)$ does not  satisfy the symmetric order condition for  
$(<_D,<_A)$. Thus there exist a bispecial   $v$  in $F(w)\cup F(\bar w)$ and letters $a$, $a'$, $b$, $b'$, such that  $a<_A a'$, $b>_D b'$ or equivalently $b<_A b'$, $avb$ and $a' vb'$ are in $F(w)\cup F(\bar w)$. Thus necessarily $a'vb'$ is in $F(w)$ and $avb$ is  in $F(\bar w)$ (or the 
opposite), thus $a'vb'$ and $b\bar v a$ are in $F(w)$ (or else  $avb$ and $b'\bar v a'$  but then  we can exchange $v$ and  $\bar v$).
By Proposition \ref{pfl}, we can build a language $L$ satisfying the symmetric order condition and containing $F(w)$. Then $v$ is  bispecial in $L.$  

Suppose first $v\neq \bar v$. By Theorem \ref{rich}  $L$ is a rich language, and by Propositions 2.3 and 2.4 of \cite{bdgz}, every complete return word of $v$
contains one occurrence of $\bar v$, and, if a word  $z$ begins with $v$, ends with $\bar v$, and has no other occurrence of $v$ or $\bar v$, then $z$ is a palindrome; and these two assertions are still true if we exchange $v$ and $\bar v$.

We show now that in $L$ there is no word $z$ beginning with  $vb'$ and ending with $av$. Indeed, if $z$ contains no other occurrence of $v$ or $\bar v$, this is impossible as $z$ is a palindrome and $b\neq b'$. Otherwise, $z$ contains $r$ occurrences of $v$, and $r$ occurrences of $\bar v$; suppose the $i$-th occurrence of $v$ is followed by $b_i$, $1\leq i\leq r$, and the $i$-th occurrence of $\bar v$ is followed by $a_i$, $1\leq i\leq r-1$; thus the $i$-the occurrence of $v$ is preceded by $a_{i-1}$, $2\leq i\leq r$,  and the $i$-th occurrence of $\bar v$ is preceded by $b_i$. We have $b_1=b<_A b'$, thus $b_1<_A b"$, and by the hypothesis and the symmetric order condition for $\bar v$ in $L$, we get $a_1<_A a'$; by the hypothesis and the order condition for $v$ we get $b_2<_Ab'$, and by induction we get $a_{r-1}<_Aa'$, $b_r<_Ab'$, hence $b_r$ cannot be $b$. Similarly  in $L$ there is no word $z$ beginning with  $\bar va'$ and ending with $av$. 

The same reasoning works when $v$ is a nonempty palindrome, using the fact that all its complete return words are palindromes, and when  $v$ is the empty word, by writing  a word $bb_1...b_rb'$, and using the order condition for $\varepsilon$ at each step. 

Thus no word in $L$ can contain both $avb$ and $b'\bar va'$, thus $w$ cannot be  in $L$, contradiction. \end{proof}

\section{postscript} During the course of submission of the present article, we discovered a second paper by F. Dolce and C.B. Hughes \cite{d2} which also fully answers  M. Lapointe's question in \cite{lap}. In fact, the main result of their paper states that each return word $u$ in the natural coding of a standard interval exchange transformation for a pair of orders $(<_D,<_A)$ clusters for  $(<_D,<_A)$  (Theorem 6.1 in ~\cite{d2}).  In the actual paper it is stated in terms of the permutation $\pi$ rather than pair of orders. The methods are  different as they rely on an extended form of branching Rauzy induction whose theory is developed in their paper. As standard interval exchange transformations form a sub-class of all generalised  interval exchange transformations, Theorem 6.1 in \cite{d2} is a consequence of our Theorem~\ref{simple}. But in fact the converse is also true: if a primitive word $u$ clusters for a pair of orders $(<_D,<_A),$  then by \cite{fz4} the language $L_u$ is generated by a standard interval exchange transformation with departure and arrival orders $(<_D,<_A)$ and $u$ is a return word to itself in the language $L_u.$

\end{document}